\documentclass[11pt]{amsart}

\usepackage{amsfonts}
\usepackage{amssymb}
\usepackage{amscd}
\input xy
\xyoption{all}

\def\ZZ{{\mathbb Z}}

\def\QQ{{\mathbb Q}}

\def\PP{{\textbf P}}
\def\OO{{\mathcal O}}
\def\R{\mathbf{R}}
\def\LL{\mathbf{L}}
\def\SS{\mathcal{S}}
\def\D{\mathbf{D}}

\def\F{\mathcal{F}}
\def\E{\mathcal{E}}
\def\G{\mathcal{G}}
\def\H{\mathcal{H}}

\def\I{\mathcal{I}}
\def\J{\mathcal{J}}

\def\Pic0{{\rm Pic}^0(X)}

\theoremstyle{plain}

\newtheorem{theorem}{Theorem}[section]

\newtheorem{proposition/example}[theorem]{Proposition/Example}
\newtheorem{proposition}[theorem]{Proposition}
\newtheorem{corollary}[theorem]{Corollary}
\newtheorem{lemma}[theorem]{Lemma}

\theoremstyle{definition}
\newtheorem{definition}[theorem]{Definition}
\newtheorem{remark}[theorem]{Remark}
\newtheorem{example}[theorem]{Example}

\newtheorem{conjecture/question}[theorem]{Conjecture/Question}

\newtheorem{remark/definition}[theorem]{Remark/Definition}
\newtheorem{definition/notation}[theorem]{Definition/Notation}

\newlength{\sectiontitlewidth}
 \theoremstyle{remark}

\begin{document}

\title{Generic Vanishing filtrations and perverse objects in derived categories of coherent sheaves}

\author[M. Popa]{Mihnea Popa}
\address{Department of Mathematics, University of Illinois at Chicago,
851 S. Morgan Street, Chicago, IL 60607, USA } \email{{\tt
mpopa@math.uic.edu}}
\thanks{The author was partially supported by the NSF grant DMS-0758253 and a Sloan Fellowship.}

\date{\today}
\maketitle

\tableofcontents

\setlength{\parskip}{.1 in}

\markboth{M. POPA} {\bf Generic Vanishing filtration on derived categories}

\section{Introduction}

This is a mostly expository paper, intended to explain a very natural relationship between two a priori 
distinct notions appearing in the literature: \emph{generic vanishing} in the context of vanishing theorems and birational geometry (\cite{gl1}, \cite{gl2}, \cite{clh}, \cite{hacon}, \cite{pp2}, \cite{pp4}), and \emph{perverse coherent sheaves} in the context of derived categories (\cite{kashiwara}, \cite{bezrukavnikov}, \cite{abe}, \cite{yz}). Criteria for checking either condition are provided and applied 
in geometric situations. A few new results, and especially new proofs, are included as well.

Let $X$ and $Y$ be noetherian schemes of finite type over a field, with $\D(X)$ and $\D(Y)$ denoting their 
bounded derived categories of coherent sheaves. If $X$ is proper and $P$ is a perfect object in 
$\D (X\times Y)$, 
we have an integral (or Fourier-Mukai-type) functor
$$\R \Phi_P : \D (X) \rightarrow \D(Y),
~~ \R \Phi(\cdot) : = \R {p_Y}_* ( p_X^*(\cdot) \overset{\LL}{\otimes}
P).$$

Based on \cite{pp2} and \cite{pp4} one can introduce a Generic Vanishing ($GV$) filtration on the derived category $\D(X)$, with respect to this functor. Denote $\dim X = d$ and $\dim Y = g$. 

\noindent
{\bf Definition.}
For every integer $m$, define the full subcategory of $\D(X)$
$$GV_m (X) := \{A \in \D(X) ~|~ {\rm~codim}~{\rm Supp} ~R^i \Phi_P A \ge i + m {\rm~for~all~} i > 0\}.$$
We say that an object in $GV_m (X)$ \emph{satisfies Generic Vanishing with index $m$} with respect to $P$.
(The definition depends of course on $P$, and should rather be $GV_m^P (X)$, but unless there is danger of 
confusion I will avoid this to simplify the notation.)

\noindent
For dimension reasons we have
$$GV_g(X) = GV_{g+1} (X) = \ldots = \R\Phi_P^{-1} (\D^{\le 0} (Y)).$$ 
On the other hand, the negative indices can go indefinitely, as for $k\ge 0$ we have that $A \in GV_{- k}(X)$ is equivalent to $A[-1] \in GV_{-k-1}(X)$, i.e.
$$ GV_{-k} (X) = GV_0 (X)  [-k], ~{\rm for~all~} k \ge 0.$$
(This is not the case for $k <0$.) In conclusion
the integral functor $\R\Phi_P$ induces a \emph{Generic Vanishing filtration} on $\D(X)$ given by 
$$\R\Phi_P^{-1} (\D^{\le 0} (Y))= \ldots = GV_g (X)  \subset \ldots \subset GV_1 (X) 
\subset GV_0 (X) =: GV(X) \subset$$
$$\subset  GV_{-1} (X) \subset \ldots \subset GV_{-d}(X) \subset GV_{-d - 1} (X) \subset \ldots$$  

\noindent
In \S2 I compare this with a cohomological support loci filtration appearing in vanishing theorems. 

By definition, the images of objects in the central component $GV (X)$ correspond via $\R\Phi_P$ to the negative component of a perverse $t$-structure on $\D(Y)$ introduced by Bezrukavnikov \cite{bezrukavnikov} (following Deligne; cf. also  \cite{abe}) and Kashiwara \cite{kashiwara},\footnote{This is one of many perverse $t$-structures, each depending on the choice of a perversity function. The main result of Bezrukavnikov and Kashiwara is to 
describe \emph{all} functions that do define such $t$-structures. In this paper the notion of perversity corresponds only to this particular $t$-structure (cf. Definition \ref{function}), which is very geometric in flavor, but only a sliver of the algebraic picture. It is an interesting question whether there exist natural geometric objects satisfying generic vanishing conditions associated to other perversity functions, especially 
since a few such functions have recently appeared in the study of Donaldson-Thomas invariants  \cite{bayer}, \cite{toda1}, \cite{toda2} or stability conditions \cite{meinhardt}.}  Grothendieck dual to the standard $t$-structure. I include in \S3 a self-contained proof of the existence of this $t$-structure, following methods in \cite{kashiwara} and \cite{pp4}, and identify the preimage of its positive part in $\D(X)$ via $\R \Phi_P$ as well. Perverse sheaves correspond via this functor to what will be called \emph{geometric $GV$ objects}, and furthermore via duality to objects satisfying the Weak Index Theorem with index $d$ ($WIT_d$) with respect to a related kernel; see Theorem \ref{wit_perverse}. (The standard example of such an object 
is the canonical bundle $\omega_X$ of  a variety $X$ with generically finite Albanese map, according to a fundamental theorem of Green-Lazarsfeld \cite{gl1}.) In case $\R\Phi_P$ is a Fourier-Mukai equivalence, this means that the subcategories of geometric $GV$ objects and $WIT_d$ objects in $\D(X)$ are hearts of natural $t$-structures.

The main emphasis is that there are explicit cohomological criteria for detecting all of the individual components of the Generic Vanishing filtration, or the abelian category of perverse coherent sheaves and the components of a refined filtration on it. It is shown in \S4 and \S5 how $GV_m(X)$ can be characterized either in terms of the vanishing of certain 
hypercohomology (for $m \le 0$), or in terms of commutative algebra (for $m >0$). 
While quite basic theoretically, these criteria have concrete geometric applications which provided the initial motivation for undertaking this study. Below is a brief explanation. 

I characterize in \S4 the negative part of the $GV$ filtration, when $Y$ is projective,  by the vanishing of suitable hypercohomology groups. This is a slight generalization of the main technical results in \cite{pp2}. These were in turn inspired by the  approach to the Green-Lazarsfeld results for the canonical bundle in Hacon's paper \cite{hacon}, where derived category techniques made their first appearance in the study of generic vanishing. The proof given here clarifies the treatment in \cite{pp2}: after applying Grothendieck Duality to the integral functor, it emphasizes a distinct local part having to do with Kashiwara's characterization of the perverse $t$-structure on $\D(Y)$, and a global one (for $Y$ projective) which is a simple cohomological characterization of perverse sheaves. 

Section 5 deals with the positive part of the filtration. It contains results that were announced (and proved for the correspondence 
between a smooth projective variety $X$ and $\Pic0$) in \cite{pp4}. One introduces
a filtration on the abelian category of perverse coherent sheaves given by a variant of Serre's condition 
$S_k$. Due to a criterion of Auslander-Bridger \cite{ab}, this corresponds via $\R\Phi_P$ to the 
$GV_m$ filtration with $m >0$, and in the case of finite homological dimension with a filtration by 
syzygy conditions as well. In this last case, one can apply the Evans-Griffith Syzygy Theorem to bound the 
rank of the dual of a non-locally free perverse sheaf. Significantly, this rank is sometimes a standard invariant of $X$, like its holomorphic Euler characteristic $\chi(\omega_X)$ in the case of the Fourier-Mukai transform of 
$\omega_X$. The homological commutative algebra needed here, as well as in \S3, is reviewed in an Appendix in \S7.

In \S6 I present geometric applications of the two main criteria characterizing components of the $GV$ filtration. The hypercohomology vanishing in Theorem \ref{negative} can often be checked in practice by reducing it to standard Kodaira-Nakano-type vanishing theorems. On projective varieties this accounts for essentially all known extensions of 
the generic vanishing theorems of \cite{gl1} (cf. \S6.1 and \S6.2). It can also be applied in connection with problems related to moduli spaces of vector bundles (cf. \S6.5) and to cohomology classes on abelian varieties (cf. \S6.4). On the other hand, the syzygy criterion in Theorem \ref{syzygy_fm} is applied in conjunction with the Evans-Griffith theorem as expalained in the previous paragraph, in situations when the cohomological support loci are known to be small, to the study of irregular varieties (cf. \S6.3).

At the moment there is an obvious disconnect between the theoretical results in \S3-5, which work in a 
general setting, and the applications in \S6 and elsewhere, which are almost all in the context of the integral functor induced by a universal line bundle on $X \times \Pic0$, for a smooth projective $X$. While this is not very restrictive from the point of view of 
birational geometry, it is natural to expect generic vanishing phenomena (or, equivalently, interesting perverse sheaves) on other moduli spaces, and also on some spaces with singularities. I give a few such examples in \S6.5, but full results in this direction are still to be discovered, the main obstruction being the current poor understanding of moduli spaces of sheaves on varieties of dimension three or higher.\footnote{More precisely, from the present perspective one needs a good understanding of the Fourier-Mukai transforms of their determinant line bundles.}
Hence overall this material has two rather distinct aspects: characterizing $t$-structures defined by support conditions in a formal setting on one hand, and applying this to the concrete geometric study of generic vanishing for the Picard variety (or other parameter spaces to a lesser extent) on the other.  Work in the two directions has been done by somewhat disjoint groups of people. I hope this note indicates that the results and methods involved often overlap or are extremely similar, and will serve as an introduction to the geometric aspects for those more algebraically inclined and viceversa.

\noindent
{\bf Acknowledgements.} 
Much of the material  that is due to the author has been worked out in articles or discussions 
with G. Pareschi, so this paper should be considered at least partially joint with him. I also thank Dima Arinkin,  David Ben-Zvi, Iustin Coand\u a, Daniel Huybrechts, Robert Lazarsfeld, Mircea Musta\c t\u a and Christian Schnell for very useful conversations. It was 
David Ben-Zvi who first pointed out \cite{bezrukavnikov} and Daniel Huybrechts who first pointed out 
\cite{kashiwara}, both guessing that our work should have a connection with these papers.  Some of the material was presented at a workshop at University of Michigan in May 2009. Special thanks are due to the organizers, especially Mircea Musta\c t\u a, for the opportunity and for an extremely pleasant week.

\section{Preliminaries and examples}

\noindent
{\bf Integral functors.}
The technical results will be proved for noetherian schemes of finite type over a field.\footnote{ 
Note that all of the results in \S3--6, except for those in \S6 involving the Evans-Griffith
syzygy theorem, hold more generally for schemes of finite type over the spectrum of a ring 
with dualizing complex. The reader will also note that many of the results hold, with 
appropriate formulations, for complex analytic spaces.} 
Given any such $X$, we denote by $\D(X)$ the bounded derived category of coherent 
sheaves on $X$. For any other $Y$ of this type, and $P$ a perfect object in $\D (X\times Y)$ 
(or more generally of finite tor-dimension over $X$ and $Y$),  if $X$ is proper we have an integral functor
$$\R \Phi_P : \D (X) \rightarrow \D(Y),
~~ \R \Phi(\cdot) : = \R {p_Y}_* ( p_X^*(\cdot) \overset{\LL}{\otimes}
P).$$
If $Y$ is also proper, we have the analogous
$$\R \Psi_P : \D(Y) \rightarrow \D (X),
~~\R \Psi(\cdot) : = \R {p_X}_* ( p_Y^*(\cdot) \overset{\LL}{\otimes} P).$$
According to standard terminology, an object $A$ in $\D(X)$ is said to satisfy $WIT_j$, i.e. the Weak Index 
Theorem with index $j$, if $R^i \Phi_P A = 0$ for all $i \neq  j$.

The dual of an object $A$ in $\D(X)$ is $A^\vee : = \R \mathcal{H}om(A, \OO_X)$. 
When $X$ is Cohen-Macaulay, we also denote $\R \Delta A :=  \R \mathcal{H}om ( A, \omega_X)$, 
where $\omega_X$ is the dualizing sheaf of $X$. When $X$ is in addition projective, 
Grothendieck duality takes the following form (cf. \cite{pp2} Lemma 2.2):

\begin{lemma}\label{gd}
Assume that $X$ is Cohen-Macaulay and projective, of dimension $d$. Then for any $A$ 
in $\D(X)$, the 
Fourier-Mukai and duality functors satisfy the exchange formula
$$(\R\Phi_P A)^\vee \cong \R{\Phi}_{P^\vee} (\R \Delta A) [d].$$
\end{lemma}

\noindent
{\bf $t$-structures.} Recall the following well-known
\begin{definition}
Let $\D$ be a triangulated category. A $t$-structure on $\D$ is a pair of full subcategories  $\big(\D^{\le 0}  , \D^{\ge 0} \big)$, with 

$\D^{\le n} := \D^{\le 0} [-n]$ 
and $\D^{\ge n} := \D^{\ge 0} [-n]$, satisfying

\noindent
(a) $\D^{\le 0} \subset \D^{\le 1}$ and $\D^{\ge 1} \subset \D^{\ge 0}$. 

\noindent
(b) ${\rm Hom}_{\D} (A, B) = 0$, for all $A \in  \D^{\le 0}$ and $B\in \D^{\ge 1}$.

\noindent
(c) For any $X\in \D$ there exists a triangle
$$A \rightarrow X\rightarrow B\rightarrow A[1]$$
with $A \in  \D^{\le 0}$ and $B\in \D^{\ge 1}$.
\end{definition}

\noindent
The category $\D^{\le 0}   \cap \D^{\ge 0}$ is called the \emph{heart} (or \emph{core}) of the $t$-structure; it is an abelian category. The main example is of course, for a scheme $X$ as above, the \emph{standard $t$-structure} on $\D(X)$ given by $\big(\D^{\le 0} (X)  , \D^{\ge 0} (X) \big)$. Its heart is ${\rm Coh} (X)$. Another important $t$-structure will appear in \S3.

\noindent
{\bf Cohomological support loci.}
Generic vanishing conditions were originally given in terms of cohomological support
loci, so it is natural to compare the definition of a $GV_m$-object $A$ in the Introduction
with the condition that the $i$-th cohomological support locus of
$A$ has codimension $\ge i +m$. For any $y\in Y$ we denote $P_y= \LL
i_{y}^* P$ in $\D(X_y)= \D(X)$, where $i_{y}: X_y = X\times \{y\}
\hookrightarrow X \times Y$ is the inclusion.

\begin{definition}
Given an object $A$ in $\D(X)$, the \emph{$i$-th cohomological support
locus} of $A$ with respect to $P$ is
$$V^i_P(A) : = \{y\in Y\>|\>  H^i(X_y, A \overset{\LL}{\otimes} P_y) \neq 0\}.$$
\end{definition}

Although ${\rm Supp}~R^i\Phi_P A$ and $V^i_P (A)$ are in general different, they carry the same numerical information in the following sense:

\begin{lemma}\label{equivalence}
For every $m \in \ZZ$, the following conditions are equivalent:\\
(1) $A$ is a $GV_m$-object with respect to $P$. \\
(2) ${\rm codim}~V^i_P(A)\geq i +m$ for all $i$.
\end{lemma}
\begin{proof}
I reproduce the proof in \cite{pp2} Lemma 3.6, for the sake of completeness.
Since by cohomology and base change 
(for the hypercohomology of bounded complexes -- cf. \cite{ega3} 7.7, especially 7.7.4, 
and Remarque 7.7.12(ii)) we have that ${\rm Supp} ~R^i \Phi_P A 
\subseteq V^i_P(A)$, it is enough to prove that (1) implies (2).
The proof is by descending induction on $i$. There certainly exists
an integer $s$ such that $H^{j}(X, A \overset{\LL}{\otimes} P_y)=0$ for any
$j>s$ and for any $y\in Y$. Then, by base change, ${\rm Supp}~R^s
\Phi_P A =V^s_P(A)$. The induction step is as follows: assume that
there is a component $\bar V$ of $V^i_P(A)$ of codimension less
than $i+m$. Since (1) holds, the generic point of $\bar V$ cannot be
contained in ${\rm Supp}~R^i \Phi_P A$ and so, again by base
change, we have that $\bar V\subset V^{i+1}_P(A)$. This implies
that ${\rm codim}~ V^{i+1}_P(A)< i + m$, which contradicts the
inductive hypothesis.
\end{proof}

\noindent
{\bf Examples.} Here are some basic examples related to the Generic Vanishing filtration. Many more examples will appear in \S6.

\noindent
(1) The pioneering result on generic vanishing is the following theorem of Green-Lazarsfeld, \cite{gl1} Theorem 1: if $X$ is a smooth projective variety and the general fiber of its Albanese map has dimension $k$, then 
$\omega_X \in GV_{-k} (X)$. This is of course with respect to the integral functor $\R\Phi_P : \D(X) \rightarrow 
\D(\Pic0)$ given by a normalized Poincar\'e bundle on $X\times \Pic0$. If $k = 0$ we will see below criteria for when 
$\omega_X \in GV_m (X) $ for some $m >0$.

\noindent
(2)  Let $C$ be a smooth projective curve of genus $g \ge 2$, and $P$ a Poincar\'e bundle on $C \times {\rm Pic}^0 (C)$.
It is a simple exercise to check that if $L$ is a line bundle of degree $d$ on $C$, then $V^i (L) = \emptyset$ for $i \ge 2$ and 
$$V^1 (L) \cong W_{2g-2-d}(C) \subset {\rm Pic}^{2g-2-d}(C),$$
the image of the Abel-Jacobi map of the $(2g-2 -d)$-th symmetric product of $C$. This has dimension $2g-2-d$ for $d \ge g-1$. 
Hence the restriction of the Generic Vanishing filtration to ${\rm Pic} (C)$ looks as follows
$$GV_g^* (C)  \subset \ldots \subset GV_1^* (C)  \subset GV_0^* (C)  \subset  GV_{-1}^* (C)$$
where $GV_m^* (C) := GV_m (C) \cap {\rm Pic}(C)$ is equal to the whole of ${\rm Pic}(C)$ for $m = -1$ and to 
${\rm Pic}^{\ge g-1 + m} (C)$ for $m \ge 0$. 

\noindent
As a side remark, note that one can similarly introduce more refined loci 
$$V^1_p (L) = \{\alpha\in {\rm Pic}^0(C) ~|~ h^1(C, L \otimes \alpha) \ge p\}.$$
As we vary $d$ and $p$, these correspond to all Brill-Noether loci $W_d^r (C)$ (see e.g. \cite{acgh} Ch.IV, \S3), which can be quite complicated.
An understanding of the pieces of the Generic Vanishing filtration on $\D(C)$, and of suitable refinements, could then be seen as a broad generalization of Brill-Noether theory to arbitrary coherent sheaves, or even objects in the derived category. The same goes for $\D(X)$ for higher dimensional $X$, where even the line bundle picture is quite mysterious.

\section{Comparison with a Bezrukavnikov-Kashiwara perverse $t$-structure}

Let $Y$ be a Cohen-Macaulay scheme of finite type over a field. Bezrukavnikov \cite{bezrukavnikov} and Kashiwara \cite{kashiwara} have introduced a $t$-structure on $\D(Y)$ which corresponds to the standard $t$-structure via taking derived duals, as part of a general procedure of defining perverse $t$-structures on the bounded derived category of coherent sheaves.\footnote{In fact Bezrukavnikov works under slightly more general, and Kashiwara under slightly more restrictive hypotheses.} The exposition here is  closer to that of \cite{kashiwara}. Explicitly, define
$${}^p\D^{\le 0} (Y) := \{ A \in \D(Y) ~|~ {\rm codim}~{\rm Supp}~\mathcal{H}^i A \ge i {\rm~ for~all~} i \ge 0\} ~{\rm and~}$$ 
$${}^p\D^{\ge 0} (Y) := \{ A \in \D(Y) ~|~ \mathcal{H}^i_Z(A) = 0, ~\forall Z\subset X {\rm~closed ~with~}
{\rm codim}~Z > i \}.$$
Here $\mathcal{H}^i_Z(\cdot)$ denotes local hypercohomology with support in $Z$. The following was proved in 
\cite{kashiwara} Proposition 4.3 in the smooth case (cf. also \cite{bezrukavnikov} Lemma 5, in a more 
general setting):

\begin{proposition}\label{kashiwara}
(1) The pair $\big( {}^p\D^{\le 0} (Y), {}^p\D^{\ge 0} (Y)\big)$ is a $t$-structure on $\D(Y)$.

\noindent
(2) We have $\R\Delta \big(\D^{\ge 0} (Y)\big) =  {}^p\D^{\le 0} (Y)$ and $\R\Delta \big(\D^{\le 0} (Y)\big) =  {}^p\D^{\ge 0} (Y)$.
\end{proposition}

Regarding (2), the proof in \cite{kashiwara} shows the inclusions from left to right. Together with the fact that the pair $\big( {}^p\D^{\le 0} (Y), {}^p\D^{\ge 0} (Y)\big)$ is a $t$-structure (cf. \cite{kashiwara} Theorem 3.5), this formally implies that they are both equalities.
To keep the exposition self-contained, I give below direct proofs for both equalities in (2) for $Y$ as above, using arguments from \cite{pp4} \S2 and \cite{kashiwara} \S4 plus some consequences of duality for local cohomology listed in Corollary \ref{support}. It implies that (1) also holds, by transfer of the standard $t$-structure via duality. The two equalities are treated separately in the following two Propositions, stated slightly more generally.

\begin{proposition}\label{cm1}
For every $k \ge 0$ we have $\R\Delta \big(\D^{\ge -k} (Y)\big) =  {}^p\D^{\le k} (Y)$, i.e. for $A \in \D(Y)$ 
$$ {\rm codim}~{\rm Supp} ~\mathcal{H}^i A \ge i -k {\rm~for~all~} i > 0 \iff R^i \Delta A = 0 {\rm ~for~all~} i  < -k.$$
\end{proposition}
\begin{proof}
Since the condition $ {\rm codim}~{\rm Supp}(\mathcal{H}^i A) \ge i -k$ is non-trivial only for $i \ge k$, by replacing $A$ by $A[k]$ we can reduce to the case $k =0$. 

Start first with $A$ such that $ {\rm codim}~{\rm Supp}(\mathcal{H}^i A) \ge i$ for all $i$. 
We have a spectral sequence 
$$E_2^{pq} : = \mathcal{E}xt^p (\mathcal{H}^q A, \omega_Y) \Rightarrow R^{p - q} \Delta A.$$
In other words, to converge to $R^i \Delta A$, we start with the $E_2$-terms  $\mathcal{E}xt^{i +j} (\H^j A, \omega_Y)$. Since  ${\rm codim}~ {\rm Supp} ~\mathcal{H}^j A \ge j$ for all $j$, by Corollary \ref{support}(a) we get $\mathcal{E}xt^{i +j} (\mathcal{H}^j A, \omega_Y) = 0$ for all $i < 0$, which by the spectral sequence implies that $R^i \Delta A = 0$ for $i < 0$.

Now start with $A$ such that $R^i \Delta A = 0$ for $i < 0$.
Since $\R \Delta$ is an involution on $\D(Y)$, we have a spectral sequence 
$$E_2^{pq} : = \mathcal{E}xt^p(R^q \Delta A , \omega_Y) \Rightarrow \mathcal{H}^{p -q} A.$$
In other words to converge to $\mathcal{H}^i A$, we start with the $E_2$-terms $\mathcal{E}xt^{i +j}( R^j \Delta A , \omega_Y)$. But Corollary \ref{support}(b) implies that for all $j$ we have 
$${\rm codim}~{\rm Supp}~\mathcal{E}xt^{i +j}( R^j \Delta A , \omega_Y) \ge i +j.$$ 
Since by hypothesis we must have $j \ge 0$, this implies that the codimension of the support of all 
$E_2$-terms is at least $i$. But then chasing through the spectral sequence this immediately 
implies that ${\rm codim}~ {\rm Supp} ~\mathcal{H}^i A \ge i$.
\end{proof}

\begin{proposition}\label{cm2}
For every $k \ge 0$ we have $\R\Delta \big(\D^{\le -k} (Y)\big) =  {}^p\D^{\ge k} (Y)$, i.e. for $A \in \D(Y)$ 
$$ \mathcal{H}^i_Z(A) = 0, {\rm~for~all~} Z\subset X {\rm~closed ~with~}
{\rm codim}~Z > i- k\iff R^i \Delta A = 0 {\rm~for~all~} i  > -k.$$
\end{proposition}
\begin{proof}
One can again reduce to the case $k=0$ by shifting by $k$. 
Let $A$ such that $R^i \Delta A = 0$ for $i>0$. Denote $B = \R  \Delta A$. For any closed $Z$, we use the derived local 
duality isomorphism
$$\R \Gamma_Z (A) \cong \R \mathcal{H}om (B, \R\Gamma_Z (\omega_Y)).$$
There is a (double) spectral sequence computing the right hand side, namely
$$ \mathcal{E}xt^p (H^i B  , \mathcal{H}^j_Z(\omega_Y)) \Rightarrow 
R^{p- i + j} \mathcal{H}om (B, \R\Gamma_Z (\omega_Y)).$$
Now $H^i B = 0$ for $i >0$, while $\mathcal{H}^j_Z (\omega_Y) = 0$ for $j < {\rm codim}~Z$ by Corollary \ref{support}(a). This implies that $\mathcal{H}^i_Z (A) = 0$ for $i < {\rm codim}~Z$ as well.

Let now $A$ such that $\mathcal{H}^j_Z (A) = 0$ for any closed $Z$ such that $i < {\rm codim}~Z$.
By \cite{kashiwara} Proposition 4.6, this is equivalent to the fact that $A$ can be represented by a bounded complex 
$F^\bullet$ of locally free $\OO_Y$-modules in non-negative degrees. But then we have a spectral sequence with $E_2$-terms $\mathcal{E}xt^{i +j} (F^j , \omega_Y)$ converging to $R^i \Delta A$. Now $j \ge 0$,  hence 
$\mathcal{E}xt^{i +j} (F^j , \omega_Y) = 0$ for $i >0$, since the $F_j$ are flat. This implies that $R^i \Delta A = 0$ for $i>0$.
\end{proof}

\begin{definition}\label{function}
Coherent perverse $t$-structures were 
defined and studied in general in \cite{bezrukavnikov} Theorem 1 and \cite{kashiwara} Theorem 5.9, by means of perversity (or support) functions. 
The particular $t$-structure in Proposition \ref{kashiwara} corresponds to the perversity function 
$$p: \{0, \ldots , g\} \rightarrow \ZZ, ~p(m) = g-m,$$ 
where $g = \dim Y$ (or equivalently $p^\prime: Y \rightarrow \ZZ$, $p^\prime (y) = \dim \OO_{Y,y}$).
It is further studied via rigid dualizing complexes in a more general context in \cite{yz},
where it is called the \emph{rigid} perverse $t$-structure.
A \emph{perverse coherent sheaf} on $Y$ is an object in the heart of the $t$-structure \\ $\big( {}^p\D^{\le 0} (Y), {}^p\D^{\ge 0} (Y)\big)$. We denote these by
$${\rm Per} (Y) =  {}^p\D^{\le 0} (Y) \cap  {}^p\D^{\ge 0} (Y).$$
By Proposition \ref{kashiwara} we have simply that ${\rm Per} (Y) = \R\Delta ({\rm Coh}(Y))$. 
\end{definition}

\noindent
{\bf Comparison via integral functors.}
Let now in addition $X$ be a Cohen-Macaulay scheme, projective over $k$, 
and let $P$ be a perfect object in $\D(X\times Y)$.
The goal is to use the above discussion after applying the functor $\R \Phi_P$ to objects in $\D(X)$. When 
$Y$ is projective, the next result is the equivalence of (1) and (3) in \cite{pp2} Theorem 3.7.
In general it was stated (in the smooth case) in \cite{pp4} Theorem 5.2 -- the proof was given there only in a special case, but the method is similar. Denote
$$Q : = P^\vee \otimes p_Y^* \omega_Y.$$

\begin{corollary}\label{gv_wit}
Let $A \in \D(X)$ and $k \ge 0$. Then 
$$A \in GV_{-k}(X) \iff \R\Phi_{Q} (\R\Delta A) \in \D^{\ge d-k} (Y).$$
\end{corollary}
\begin{proof}
First, it is clear by definition that for $k\ge 0$ we have 
$$A\in GV_{-k} (X) \iff \R\Phi_P A \in {}^p\D^{\le k} (Y).$$ 
This result is then simply a reinterpretation of the equivalence in Proposition \ref{cm1} applied to $\R\Phi_P A$, 
via duality. Indeed, 
by Lemma \ref{gd} we have that $\R\Phi_{Q} (\R\Delta A) \in \D^{\ge d-k} (Y)$ if and only if 
$(\R \Phi_{Q^\vee} A)^\vee \in \D^{\ge -k} (Y)$. But the projection formula implies 
$$(\R \Phi_{Q^\vee} A)^\vee \cong \R \Delta (\R\Phi_P A).$$
\end{proof}

\begin{example}[{\bf The Fourier transform of the canonical bundle}]\label{gl_perverse}
In the situation of Example (1) in \S2, we saw that $\omega_X \in GV_{-k} (X)$ (with respect to a Poincar\'e bundle 
on $X\times \Pic0$). On the other hand, by dimension reasons we have $R^i \Phi_{P^\vee} \OO_X = 0$ for $i > d$. 
This means that when $X$ is of maximal Albanese dimension (i.e. $k = 0$), $\omega_X \in GV (X)$ and 
equivalently by Corollary \ref{gv_wit}, $\R\Phi_{P^\vee} \OO_X$ is a sheaf supported in degree 
$d$ (see Theorem \ref{gv_main} below for a  more general statement). 
In other words, $\R\Phi_P \omega_X$ is a perverse sheaf on $\Pic0$. 
\end{example}

In concrete geometric applications it is often the case as in the Example above that $P$ is a sheaf, flat over $Y$, and 
$\R\Delta A \in \D^{\le 0} (X)$ (usually a sheaf). By basic properties of push-forwards we then have $R^i \Phi_Q (\R \Delta A) = 0$ for all $i >d$, as $d$ is the dimension of the fibers of the second projection. The $GV$ such objects deserve a special name, as they are intimately related to perverse sheaves in $\D(Y)$ via $\R\Phi_P$. 

\begin{definition}\label{geom_gv}
A \emph{geometric $GV$-object} on $X$ is an object $A\in GV(X)$ such that 
$$R^i \Phi_Q (\R \Delta A) = 0  {\rm ~for~ all~} i >d.$$
\end{definition}

The conclusion of the discussion in this section can be summarized in the following result. With respect to \cite{pp2} I am also adding the new observation that under reasonable assumptions (for example when $X$ and $Y$ are smooth and $P$ lives in non-negative degrees) the $GV$-objects in non-negative degrees are always geometric.\footnote{I thank Ch. Schnell for help with this issue.} The picture will be expanded in the next two sections.

\begin{theorem}\label{wit_perverse}
Let $X$, $Y$ and $P$ be as above. For $A \in \D(X)$, the following are equivalent:

\noindent
(1) $A$ is a geometric $GV$-object.

\noindent
(2) $\R \Phi_P A$ is a perverse coherent sheaf.

\noindent
(3) $\R \Delta A$ satisfies $WIT_d$ with respect to $Q = P^\vee \otimes p_Y^* \omega_Y$.

\noindent
Assuming in addition that $P \otimes p_X^* \omega_X \in {}^p\D^{\ge 0} (X \times Y)$ 
(i.e. it can be represented by a complex of locally free sheaves in non-negative degrees),
$$\R\Phi_P \big( GV (X) \cap \D^{\ge 0} (X) \big) \subset  {\rm Per}(Y)$$
i.e. $GV$-objects in non-negative degrees are geometric.
\end{theorem}
\begin{proof}
By definition we have $\R\Phi_P^{-1} \big({}^p\D^{\le 0} (Y) \big) = GV(X)$. On the other hand, by Corollary \ref{gv_wit}, this is equivalent to $R^i \Phi_{Q} (\R\Delta A) = 0$  for all $i < d$. For the other half, note that the 
identity in Proposition \ref{cm2} applied to $\R\Phi_P A$ implies after shifting  by $k$ that 
$$\R\Phi_P^{-1} \big({}^p\D^{\ge 0} (Y) \big) = \{A \in \D(X) ~|~  
\R\Delta(\R\Phi_P A)  \in \D^{\le 0} (Y)\}.$$ 
But note that the formula at the end of the proof of Corollary \ref{gv_wit}, combined with Lemma \ref{gd}, 
implies
$$ \R\Phi_Q (\R\Delta A)  \cong \R\Delta (\R \Phi_P A)[d]$$
which gives 
$$\R\Phi_P^{-1} \big({}^p\D^{\ge 0} (Y) \big) = \{A \in \D(X) ~|~  
R^i \Phi_{Q} (\R\Delta A) = 0  {\rm ~for ~} i > d\}.$$ 
Finally, the last statement in the Theorem uses part (2) of Proposition \ref{more} below.
\end{proof}

\begin{proposition}\label{more}
(1) If $P \in \D^{\le 0} (X \times Y)$, then  $\R\Phi_P \big( {}^p\D^{\le 0} (X) \big) \subset \D^{\le d}(Y)$.

\noindent
(2) If $P\otimes p_X^* \omega_X \in {}^p\D^{\ge 0} (X \times Y)$, then $\R\Phi_P \big( \D^{\ge 0} (X) \big) \subset {}^p\D^{\ge 0}(Y)$.
\end{proposition}
\begin{proof}
(1) Let $A \in  {}^p\D^{\le 0} (X)$. We have a spectral sequence
$$E_2^{pq} : = R^p \Phi_P (\mathcal{H}^q A) \Rightarrow R^{p+q} \Phi_P  A.$$
By definition we have ${\rm dim}~{\rm Supp}~\mathcal{H}^q A \le d - q$ for all $q$.  
Since $P$ is supported only in negative degrees, we get by basic properties of push-forwards that
$R^p \Phi_P (\mathcal{H}^q A) = 0$ for all $p > d -q$, which implies the result.

\noindent
(2) Let $A \in \D^{\ge 0} (X)$. By Proposition \ref{cm1} we have $\R\Delta A \in {}^p\D^{\le 0} (X)$. 
Since $\omega_Y \overset{\LL}\otimes \omega_Y^\vee \cong \OO_Y$, by Proposition \ref{cm2}  we have $Q \in \D^{\le 0} (X \times Y)$. By (1) we then have $\R\Phi_Q (\R\Delta A) \in \D^{\le d} (Y)$. 
But we noted in the proof of the Theorem that
$\R\Phi_Q (\R\Delta A)  \cong \R\Delta (\R \Phi_P A)[d]$.
We deduce $\R\Delta (\R \Phi_P A) \in \D^{\le 0} (Y)$, which by Proposition \ref{cm2} is equivalent to 
$\R \Phi_P A \in {}^p\D^{\ge 0} (Y)$.
\end{proof}

\begin{remark}[{\bf Gorenstein asumption}]
When $Y$ is Gorenstein, the right hand side condition in Corollary \ref{gv_wit}, the definition of  a 
geometric $GV$-object, and condition (3) in Theorem \ref{wit_perverse} can all be phrased 
in terms of the simpler looking $\R\Phi_{P^\vee} (\R \Delta A)$ 
instead of $\R\Phi_Q (\R \Delta A)$. Since this will always be the case in applications, we will 
use it freely in what follows.
\end{remark}

\begin{example}
With respect to the last statement in Theorem \ref{wit_perverse}, it is not the case that the preimage in $\D(X)$ of the category of perverse sheaves is equal to 
$GV (X) \cap \D^{\ge 0} (X)$. This does not happen even in the most favorable case: let $X$ be an abelian variety, $Y = \widehat{X}$ its dual, and $P$ a Poincar\'e bundle on $X\times \widehat{X}$. Consider 
$L$ a nondegenerate line bundle on $X$ (i.e. with $\chi(L) \neq 0$) satisfying $WIT_1$, so that $L^{-1}$ satisfies $WIT_{g-1}$. 
Consider $A = L[1] \not\in \D(X)^{\ge 0}$. Then $\R\Phi_P A$ is a vector bundle supported in degree $0$, 
so $A \in GV(X)$. On the other hand $\R \Phi_{P^\vee} (\R \Delta A) \cong R^{g-1} \Phi_{P^\vee} L^{-1} [-g]$ 
is supported in degree $g$, so $\R\Phi_P A$ is perverse. 
\end{example}

A consequence of Proposition \ref{kashiwara} and its avatar Theorem \ref{wit_perverse} is extra structure on the category of geometric objects satisfying Generic Vanishing with respect to Fourier-Mukai equivalences.

\begin{corollary}
If $\R\Phi_P$ is an equivalence, the geometric $GV$-objects with respect to $P$(or dually the $WIT_d$-objects with respect to $Q$) form the heart of a $t$-structure on $\D(X)$, whose negative half consists of the category of $GV$-objects.
\end{corollary}

\noindent
It is tempting to wonder whether this still holds under weaker assumptions on $\R\Phi_P$.

\section{Hypercohomology vanishing characterization of ${\rm Per}(Y)$ and of $GV_m (X)$ 
with $m \le 0$}

The main result of this section is a characterization of the negative part of the Generic Vanishing 
filtration. It is a slightly more general  version of \cite{pp2} Theorem 3.7 and Theorem A,  in turn following 
Hacon's approach to generic vanishing for the canonical bundle in \cite{hacon}. The proof given below 
follows a strategy different from that in the papers above, emphasizing the connection with Proposition \ref{kashiwara} and the fact that the projectivity of $Y$ is needed only for a simple cohomological characterization of perverse sheaves. By a \emph{sufficiently positive} ample line bundle I always
mean a sufficiently high power of an ample line bundle, and use the notation $L \gg 0$.

\begin{theorem}\label{negative}
Let $X$ and $Y$ be Cohen-Macaulay schemes of finite type over $k$, of dimension $d$ and $g$ respectively with $X$ projective. Let $P\in \D(X\times Y)$ be a perfect object. For $A \in \D(X)$ and $k \ge 0$, the following are equivalent:

\noindent
(1) $A \in GV_{-k} (X)$.

\noindent
(2) ${\rm codim}_Y V^i_P (A) \ge i - k$ for all $i$.

\noindent 
(3) $\R \Phi_{Q} (\R \Delta A) \in \D^{\ge d-k}(Y)$.

If in addition $Y$ is also projective, they are equivalent to the hypercohomology vanishing

\noindent
(4) $H^i(X, A \overset{\LL}{\otimes} \R\Psi_{P[g]} (L^{-1}))=0$
for all $i > k$ and any $L \gg 0$ on $Y$.
\end{theorem}
\begin{proof}
The equivalence between (1) and (2) is Lemma \ref{equivalence}.
The equivalence between (1) and (3) is Corollary \ref{gv_wit}, based on Proposition \ref{cm1}. 
So we only need to prove the equivalence between (1) and (4) when $Y$ is projective. We have by definition that 
$A\in GV_{-k} (X)$ if and only if $(\R \Phi_P A)[k] \in {}^p\D^{\le 0} (Y)$. It is a standard consequence of the derived Projection Formula and Leray isomorphism that 
$$H^j (Y,\R\Phi_P A \otimes L^{-1}) \cong H^j (X, A \overset{\LL}{\otimes} \R \Psi_P (L^{-1})) \cong 
H^{j-g}(X, A \overset{\LL}{\otimes} \R\Psi_{P[g]} (L^{-1}))$$
(cf. \cite{pp2} Lemma 2.1). Hence all we need to show is the vanishing
$$H^j (Y,\R\Phi_P A \otimes L^{-1}) = 0 {\rm~for~ all~} j > g +k {\rm~and~any~} L \gg 0 {\rm~on~} Y,$$
which is a consequence of Corollary \ref{geom_gv},  and of Lemma \ref{per_coh} which describes perverse sheaves cohomologically via a standard Serre Vanishing argument.
\end{proof}

\begin{lemma}\label{per_coh}
On a Cohen-Macaulay scheme $Y$, projective over $k$ and of dimension $g$
$${\rm Per}(Y) = \{ A \in \D(Y) ~|~ H^j (Y, A \otimes L^{-1}) = 0 ~{\rm for~all~} j \neq g {\rm~and~any~} L \gg 0 {\rm~on~} Y \}.$$
More precisely:
\begin{itemize}
\item $A \in {}^p\D^{\le 0}(Y)$ if and only if 
$H^j (Y, A \otimes L^{-1}) = 0$ for all $j >g$ and any $L \gg 0$ on $Y$
\item $A \in {}^p\D^{\ge 0}(Y)$ if and only if 
$H^j (Y, A \otimes L^{-1}) = 0$ for all $j < g$ and any $L \gg 0$ on $Y$.
\end{itemize}
\end{lemma}
\begin{proof} 
It is enough to prove the two statements at the end. Note first  that $A \in {}^p\D^{\le 0}(Y)$ is equivalent to $\R \Delta A \in 
\D^{\ge 0} (Y)$ by Proposition \ref{cm1}. Now by Grothendieck-Serre duality we have 
$$H^j (Y, A \otimes L^{-1}) \cong H^{g-j} (Y, \R\Delta A \otimes L)^*.$$
We have a spectral sequence 
$$E_2^{pq} : = H^p (Y, R^q \Delta A \otimes L) \Rightarrow H^{p + q} (Y, \R\Delta A \otimes L) .$$
By Serre Vanishing, for $L \gg 0$ we have $E_2^{pq} = 0$ for all $p >0$ and all $q$, so the spectral 
sequence degenerates and 
$$H^{g-j} (Y, \R\Delta A \otimes L) \cong H^0 (Y, R^{g-j} \Delta A\otimes L).$$
Again by Serre's theorem, this is nonzero for $L \gg 0$ if and only if $R^{g-j} \Delta A \neq 0$, 
hence the assertion. The proof of the assertion for $A \in {}^p\D^{\ge 0}(Y)$ is completely analogous, using Proposition \ref{cm2}. 
\end{proof}

\noindent
The same argument shows that $\R \Phi_Q (\R \Delta A) \in \D^{\le d-k}(Y)$ (i.e. $A$ is geometric) 
if and only if 
$$H^i(X, A \overset{\LL}{\otimes} \R\Psi_{P[g]} (L^{-1}))=0$$
for all $i < k$ and all $L$ sufficiently positive on $Y$. In particular one can check the perversity of the Fourier-Mukai transform via vanishing on $X$, which is crucial for the applications in \S7.

\begin{corollary}
For $A \in \D(X)$ we have 
$$\R\Phi_P A \in {\rm Per}(Y) \iff H^i(X, A \overset{\LL}{\otimes} \R\Psi_{P[g]} (L^{-1}))=0, {\rm~for~all~} i \neq 0 {\rm~and ~all~} L \gg 0.$$
\end{corollary}

\begin{remark}
(1) The path pursued in \cite{hacon} and \cite{pp2} is to prove instead the equivalence between (3) and (4) in Theorem \ref{negative}, still based on Serre Vanishing. This equivalence is the natural extension to integral functors of a basic degeneration of the Leray spectral sequence used in the proof of Grauert-Riemenschneider-type theorems. Concretely, 
let $f: X \rightarrow Y$ be a morphism of smooth projective varieties, and consider $P : = \OO_{\Gamma}$ as a sheaf on $X\times Y$, where $\Gamma \subset X\times Y$ is the graph
of $f$. Hence $P$ induces the integral functor $\R\Phi_P = \R
f_*$, and $\R\Psi_P$ is the adjoint $\LL f^*$. Consider $A$ and $B = \R \Delta A$ objects in 
$\D(X)$. A routine calculation shows that the equivalence of (3) and (4) in Theorem
\ref{negative} applied to $A$ is the same as the well-known statement (individually for each $i$):
$$R^i f_* B = 0 \iff H^i (X, B\otimes f^* L) = 0 {\rm~for~all}~L \gg 0 {\rm ~on~}Y.$$
For instance, say $B = \omega_X$ and $f$ is generically finite. Then for any ample $L$ on $Y$, $f^* L$ is 
big and nef, so $H^i (X, \omega_X \otimes f^*L) = 0$ by Kawamata-Viehweg vanishing.
We get that $R^i f_* \omega_X  = 0 ~ {\rm~for~all~} i > 0$, 
which is of course Grauert-Riemenschneider vanishing (in the projective case).

\noindent
(2) It is interesting to note that the vanishing condition in (4) is of a different nature from standard vanishing theorems.
For instance, when $X$ and $Y$ are dual abelian varieties 
and $\R\Phi_P$ is the standard Fourier-Mukai functor, Mukai showed that $\phi_L^* \R\Psi_{P[g]} (L^{-1})\cong \oplus L$, where 
$\phi_L: Y\rightarrow X$ is the standard isogeny associated to $L$ (see the proof of Theorem \ref{gv_main} below for 
more details). This suggests that, at least when $P$ is a locally free sheaf,  $\R\Psi_{P[g]} (L^{-1})$ should be interpreted as a positive vector bundle, but which is less and less so as $L$ becomes more positive.

\end{remark}

\section{Commutative algebra filtration on ${\rm Per}(Y)$, describing $GV_m (X)$ with $m >0$}

This section is concerned with a characterization of the positive part of the Generic Vanishing 
filtration, extending the results in \cite{pp4} \S3. 
Note to begin with that the subcategories $GV_m(X)$ with $m >0$ are not obtained simply 
by shifting $GV_0 (X)$, as in the negative case. The main result, Theorem \ref{syzygy_fm}, 
is at this stage merely a matter of notation and of navigating through results in other sections and in the 
Appendix.

Let $Y$ be a Cohen-Macaulay scheme of finite type over a field, not necessarily projective. The characterization of perverse coherent sheaves in Proposition \ref{kashiwara} gives
$${\rm Per}(Y) = \{ A \in \D(Y) ~|~ \R \Delta A {\rm~is~a ~sheaf~in ~degree~} 0\}.$$
We can naturally consider a filtration on ${\rm Per}(Y)$ according to the singularities of the sheaf $\R\Delta A$.  To this end, define for $m\ge 0$ (see Definition \ref{serre}):
$${\rm Per}_m (Y) := \{ A \in {\rm Per} (Y) ~|~ \R\Delta A {\rm~satisfies~} S_m^\prime \}.$$
We get a filtration
$$\{{\rm locally~free~sheaves}\} \subset {\rm Per}_{\infty} (Y) =  \ldots =  {\rm Per}_g (Y) \subset\ldots \subset {\rm Per}_1 (Y)\subset {\rm Per}_0 (Y) = {\rm Per} (Y)$$

If we restrict to the subcategory ${\rm Per}^{{\rm fhd}} (Y) \subset {\rm Per}(Y)$ consisting of objects such that $\R \Delta A$ has finite homological dimension (so all of ${\rm Per}(Y)$ if $Y$ is smooth), then we have
$$\{{\rm maximal~Cohen-Macaulay~sheaves}\} = {\rm Per}^{{\rm fhd}}_{\infty} (Y) = \ldots 
= {\rm Per}^{{\rm fhd}}_g (Y).$$
(Note that in general $\F$ is a maximal Cohen-Macaulay  sheaf if and only if $\R\Delta\F$ is 
a maximal Cohen-Macaulay sheaf; see \cite{bh} Theorem 3.3.10.) 
On the other hand, by the Auslander-Bridger criterion in Proposition \ref{cm_fundamental}
$${\rm Per}^{{\rm fhd}} _m (Y) = \{ A \in {\rm Per}^{{\rm fhd}} (Y) ~|~ \R\Delta A {\rm~is~an~}m{\rm -th~syzygy~sheaf}\}.$$ 
Note that this last thing holds for the entire ${\rm Per}_m (Y)$ if $Y$ is Gorenstein in codimension less than or equal to one.
The Evans-Griffith Syzygy Theorem \ref{syzygy_theorem} can be rephrased as follows.

\begin{corollary}\label{perverse_rank}
Let $m> 0$ be an integer, and let $A$ be a perverse sheaf in ${\rm Per}^{{\rm fhd}}_m (Y)$ which is not a locally free sheaf. 
Then ${\rm rank} (\R\Delta A) \ge m$.
\end{corollary}

\noindent
The equivalence of (b) and (c) in Proposition \ref{cm_fundamental} gives a characterization of objects in ${\rm Per}_m (Y)$.

\begin{lemma}\label{higher_perversity}
For $A\in {\rm Per}(Y)$ and $m > 0$
$$A \in {\rm Per}_m (Y) \iff {\rm codim}~{\rm Supp} ~\mathcal{H}^i A \ge i + m {\rm~for~all~} i > 0.$$
\end{lemma}

This last condition corresponds to the $GV_m$-piece of the Generic Vanishing filtration, in the 
case of an integral functor $\R \Phi_P: \D(X) \rightarrow \D(Y)$. Recall that if $A$ is a geometric $GV$-object in $\D(X)$, then $\R \Phi_P A$ is perverse in $\D(Y)$, so $\widehat{\R\Delta A} := \R \Phi_Q (\R\Delta A)[d]$ is a sheaf.

\begin{theorem}\label{syzygy_fm}
Let $X$ and $Y$ be Cohen-Macaulay schemes of finite type over $k$, with $X$ projective. 
Fix a kernel $P \in \D(X\times Y)$, and let $A$ be a geometric $GV$-object in $\D(X)$, with respect to $P$. 
Let $m > 0$ be an integer. Then the following are equivalent:

\noindent
(1) $A \in GV_m (X)$.

\noindent
(2) $\R\Phi_P A \in {\rm Per}_m (Y)$.

\noindent
(3) $\widehat{\R\Delta A}$ satisfies $S_m^\prime$.

If these conditions are satisfied and in addition  $\R\Phi_P A \in {\rm Per}^{{\rm fhd}}_m (Y)$ or $Y$ is Gorenstein in 
codimension less than or equal to one, then they are also equivalent to 

\noindent
(4) $\widehat{\R\Delta A}$ is an $m$-th syzygy sheaf.
\end{theorem}
\begin{proof} 
The equivalence of (1) and (2) is the content of Lemma \ref{higher_perversity}. The equivalence
of (2) and (3) follows by definition and Lemma \ref{gd}. In the Gorenstein or finite homological dimension case, the 
equivalence of (3) and (4) is the Auslander-Bridger criterion quoted in Proposition \ref{cm_fundamental}.
\end{proof}

\begin{definition}[{\bf Generic Vanishing Index}]
Let $A$ be an object in $\D(X)$. The \emph{Generic Vanishing index} of $A$ (with respect to $P$) is the integer
$$gv (A) : = \underset{i>0}{\rm min}~ \{ {\rm
codim}~ {\rm Supp}~R^i \Phi_P A -i\} = \underset{i>0}{\rm min}~ \{ {\rm
codim}~ V^i_P (A) -i\}.$$ 
(The last equality holds due to Lemma \ref{equivalence}.)
If ${\rm Supp}~R^i \Phi_P A =\emptyset$ for all $i>0$, 
we declare $gv(A)=\infty$. By definition $A \in GV_m (X)$ if and only if $gv (A) \ge m$.
\end{definition}

\noindent
Theorem \ref{syzygy_fm} and Corollary \ref{perverse_rank} imply then the following useful

\begin{corollary}\label{inequality}
If $A$ is a geometric $GV$-object with $gv(A) < \infty$, then 
${\rm rank}(\widehat{\R\Delta A}) \ge gv (A)$.
\end{corollary}

\section{Geometric applications}

The characterizations of $GV_m$-objects (or of perverse objects and the syzygy filtration on them) given in \S4 and \S5 can often be  checked in practice. From a derived category point of view, one obtains nontrivial concrete examples of perverse coherent sheaves.
This produces a number of different geometric applications, some described in what  follows. The general literature on applications of generic vanishing theorems to birational geometry is very extensive, a small sampling being given by \cite{gl1}, \cite{el}, \cite{ch1}, \cite{ch2}, \cite{ch3}, \cite{hp}, \cite{pp2}, \cite{pp4}.

\noindent
{\bf 6.1. Generic Vanishing Theorems.}
The derived category approach to generic vanishing theorems was pioneered by Hacon \cite{hacon}. The work described here  is mostly taken from \cite{pp2}, and grew out of trying to extend Hacon's approach and the Green-Lazarsfeld results \cite{gl1}. 

In what follows let $X$ be a smooth projective complex variety of dimension $d$, 
with Albanese map $a: X \rightarrow A$. 
Let $P$ be a Poincar\'e line bundle on $X\times \Pic0$, and consider as usual
$$\R\Phi_P: \D(X) \rightarrow \D(\Pic0).$$ 
Every $GV$ condition will be considered with respect to this functor.
For a $\QQ$-divisor $L$ on $X$, we define $\kappa_L$ to be $\kappa(L_{|F})$, the Iitaka dimension along the generic fiber $F$ of $a$, if $\kappa(L) \ge 0$, and $0$ if $\kappa(L) = -\infty$. 

\begin{theorem}[\cite{pp2} Theorem B]\label{gv_main}
Let $L$ be a line bundle and $D$ an effective $\QQ$-divisor on $X$ such that $L-D$ is nef. If the dimension 
of the Albanese image $a(X)$ is $d-k$, then $\omega_X\otimes L \otimes \J(D)$ belongs to 
$GV_{-(k-\kappa_{(L-D)})}(X)$, where $\J (D)$ is the multiplier ideal sheaf associated to $D$.\footnote{Multiplier ideals 
are treated for example in \cite{positivity} Ch.9. Here are some easily understood examples: when $D$ is integral one has 
$\J(D) = \OO_X(-D)$, while when $D$ is a simple normal crossings $\QQ$-divisor one has $\J(D) = \OO_X (- [D])$.}  
In particular, if $L$ is a nef line bundle, then $\omega_X\otimes L$ belongs to $GV_{-(k-\kappa_L)}(X)$.
\end{theorem}

\noindent
The simplest instance of this (explaning also the terminology ``generic vanishing") is the following:

\begin{corollary}\label{nef}
Let $X$ be a smooth projective variety, and $L$ a nef line bundle on
$X$. Assume that either one of the following holds:
\begin{enumerate}
\item $X$ is of maximal Albanese dimension (i.e. $k=0$).
\item $\kappa(L)\ge 0$ and $L_{|F}$ is big, where $F$ is the generic fiber of $a$.
\end{enumerate}
Then  $\omega_X \otimes L$ belongs to $GV(X)$. In particular 
$$H^i (X, \omega_X \otimes L \otimes P) = 0 {\rm~for~all~} i > 0 {\rm~and~} P\in \Pic0 {\rm ~general}$$ 
and consequently $\chi(\omega_X \otimes L) \ge 0$.
\end{corollary}

\begin{corollary}\label{nef_perverse}
Under the hypotheses of Corollary \ref{nef}, $\R \Phi_P (\omega_X\otimes L)$ is a
perverse coherent sheaf on $\Pic0$.
\end{corollary}
\noindent
A more precise statement for $\omega_X$ is given in Corollary \ref{canonical_perverse}.

I only sketch the proof of Corollary \ref{nef}, under the hypothesis (1), as an example of the use of Theorem \ref{negative}.
This contains all the key ideas needed for Theorem \ref{gv_main}, the rest involving a standard extension to multiplier ideal sheaves, the additivity of the Iitaka dimension, and extensions of the Kawamata-Viehweg vanishing theorem 
(for full details cf. \cite{pp2} \S5). The main idea of the proof is due to Hacon \cite{hacon}, with a refinement from \cite{pp2} that allows for bypassing Hodge-theoretic results, hence the extension to twists by arbitrary nef line bundles.

Let $L$ be a nef line bundle on $X$. It is enough to show
that $\omega_X\otimes L$ satisfies condition (4) in Theorem \ref{negative}.
Let $M$ be ample line bundle on $\widehat A \cong \Pic0$, and assume for
simplicity that it is symmetric, i.e. $(-1_{\widehat A})^* M \cong
M$. We consider the two different Fourier transforms $\R \mathcal{S}
M=R^0 \mathcal{S} M = {p_A}_* (p_{\widehat{A}}^* M \otimes \mathcal{P})$ (on A), where 
$\mathcal{P}$ is a Poincar\'e bundle on $A\times \widehat{A}$ so that $P \cong (a\times {\rm id}_{\widehat{A}})^* 
\mathcal{P}$, and ${\bf R}\Psi_{P[g]} (M^{-1}) =R^g \Psi_P
(M^{-1}) =: \widehat{M^{-1}}$ (on X). These are both locally free
sheaves. One can check with a little care that
$$ \widehat{M} \cong  a^* R^g \SS (M^{-1}) \cong a^*(R^0 \mathcal{S} M)^\vee .$$
On the other hand, by \cite{mukai} 3.11, the vector bundle $R^0
\mathcal{S} M$ has the property:
$$\phi_M^*(R^0 \mathcal{S} M) \cong H^0(M)\otimes M^{-1}.$$
Here $\phi_M: \widehat A \rightarrow A$ is the standard isogeny
induced by $M$. We consider then the fiber product  $X^\prime
:=X\times_{A} \widehat A$ induced by $a$ and $\phi_M$:
$$\xymatrix{
X^\prime   \ar[r]^{\psi} \ar[d]_{b} & X
\ar[d]^{a} \\
\widehat{A} \ar[r]^{\phi_M} & A } $$ It follows that
\begin{equation}\label{square}
\psi^* \widehat{M^{-1}} \cong \psi^* a^* (R^0 \mathcal{S} M)^\vee
\cong b^*(H^0(M)\otimes M) \cong H^0(M)\otimes b^*M.
\end{equation}
What we want is the vanishing 
$$H^i (X, \omega_X \otimes L \otimes \widehat{M^{-1}})= 0 {\rm ~for~all~} i > 0.$$
(Note that this will work for \emph{any} ample line bundle $M$, the condition 
$M\gg 0$ being required only for the equivalence in Theorem \ref{negative} to hold.)
Since $\psi$, like $\phi_M$, is \' etale, it is enough to prove this after
pull-back to $X^\prime$, so for $H^i (X^\prime, \omega_{X^\prime}
\otimes \psi^* L \otimes  \psi^* \widehat{M^{-1}})$. But by
$(\ref{square})$ we see that this amounts to the vanishing
$$H^i (X^\prime, \omega_{X^\prime} \otimes \psi^*L \otimes b^*{M}) = 0 {\rm ~for~all~} i > 0.$$
Now $\psi^*L $ is nef and $b^*M$ is big and nef (as the pull-back of an ample line bundle by 
a generically finite map), so this follows from Kawamata-Viehweg vanishing.

A completely similar approach, replacing at the end Kawamata-Viehweg by other standard
vanishing theorems, proves the following results for higher direct images of canonical bundles and 
for bundles of holomorphic forms. 

\begin{theorem}\label{direct}
Let $f:Y \rightarrow X$ be a morphism, with $X$, $Y$
smooth projective varieties. Let $L$ be a nef line bundle on $f(Y)$ (reduced image of $f$). If the dimension of $f(Y)$ is $d$ and that of its image via the Albanese map of $X$ is 
$d-k$, then $R^j f_* \omega_Y \otimes L$ is a $GV_{-(k - \kappa_L)}$-sheaf on $X$ for
any $j$.
\end{theorem}

\begin{theorem}\label{nakano_1}
Let $X$ be a smooth projective variety, with Albanese image of
dimension $d-k$. Denote by $f$ the \emph{maximal} dimension of a
fiber of $a$, and consider $l:= {\rm max}\{k, f-1\}$. Then:\\
(1) $\Omega_X^j$ belongs to $GV_{-(d-j+l)}(X)$ for all $j$.\\
(2)  ${\rm codim}_{\hat A} V^i (\Omega_X^j) \ge {\rm max}\{i + j - d
-l, d - i - j -l\}$, for all
 $i$ and all $j$.\\
(3) If $L$ is a nef line bundle on $X$ and $a$ is \emph{finite},
then $\Omega^j_X\otimes L$ belongs to $GV_{-(d - j)}(X)$ for  all $j$.
\end{theorem}

Theorem \ref{direct} is again due to Hacon \cite{hacon} in the case $L = \OO_X$. Another statement for $\Omega^j_X$, that I do not know how to cover with these methods, can be found in \cite{gl1}. For the first statement one uses Koll\'ar's vanishing and torsion-freeness theorems for higher direct images of 
canonical bundles, while for the second Nakano and Bogomolov-Sommese vanishing. Details and further results along these lines can be found in \cite{pp2} \S5. Note that Theorem \ref{direct} produces other natural examples of perverse sheaves on Picard varieties.

\begin{corollary}
Let $X$ be a smooth projective variety and $a: X\rightarrow A$ its Albanese map, and let $L$ be any 
nef line bundle on $a(X)$. Then $\R\Phi_P (R^i a_* \omega_X \otimes L)$ are perverse sheaves on $\Pic0$ for all $i$, where $\R\Phi_P$ is the standard Fourier-Mukai functor on $\D(A)$.
\end{corollary}

\begin{remark}
All of these results, suitably interpreted,  hold more generally for morphisms $X\rightarrow A$, with 
$A$ an abelian variety and $\widehat{A}$ replacing $\Pic0$. This allows in particular for a more 
general statement of Theorem \ref{direct}, where $X$ can be singular. Also, it is interesting to note that the methods presented here give an algebraic proof of generic vanishing in characteristic $0$ (as well as in positive characteristic in some cases), 
which had traditionally been approached via transcendental methods. 
(For details cf. \cite{pp2}, end of \S5 and Remark 6.7.) 
\end{remark}

The most striking applications of generic vanishing theorems include the following results of Ein-Lazarsfeld \cite{el} on singularities of theta divisors, and of Chen-Hacon \cite{ch1} on a conjecture of Koll\'ar on characterizing abelian varieties (cf. also \cite{pareschi2} for a proof using the interpretation 
in \S6.2 as well).

\begin{theorem}[\cite{el} Theorem 1]
Let $(A, \Theta)$ be an indecomposable principally polarized abelian variety. Then $\Theta$ is normal
and has rational singularities.
\end{theorem}

\begin{theorem}[\cite{ch1} Theorem 3.2]
Let $X$ be a smooth projective variety with $h^0 (X, \omega_X) = h^0 (X, \omega_X^{\otimes 2}) = 1$ and $h^1(X, \OO_X) = \dim X$. Then  $X$ is birational to an abelian variety.
\end{theorem}

\noindent
{\bf 6.2. Vanishing of higher direct images.}
The equivalence between (1) and (3) in Theorem \ref{negative} says that once the (geometric) $GV$-condition is established for an object $A$,  the corresponding integral transform of the Grothendieck dual object is up to shift a sheaf. Again in the setting of $\R\Phi_P: \D(X) \rightarrow \D(\Pic0)$,  here is the main instance of this:

\begin{corollary}\label{dual}
Let $X$ be a smooth projective variety of dimension $d$ and Albanese dimension $d -k$. Then 
$$R^i \Phi_{P^\vee} \OO_X (\cong R^i {p_2}_* ( P^\vee) ) = 0 {\rm ~for~} i \not \in [d-k, d].$$
In particular, if $X$ is of maximal Albanese dimension, then $\OO_X$ satisfies $WIT_d$, so that its Fourier-Mukai transform is a sheaf 
$$\widehat{\OO_X} \cong \R\Phi_{P^\vee} \OO_X [d].$$ 
\end{corollary}

Although with the methods presented in this paper this statement is an immediate consequence of generic vanishing, its history follows somewhat the opposite direction: a more precise version was first proved when $X = A$ is an abelian variety by Mumford \cite{mumford} \S13. It was then 
conjectured to be true in general by Green-Lazarsfeld \cite{gl2} Problem 6.2. This was proved by Hacon \cite{hacon} and Pareschi \cite{pareschi1}, both showing that it further implies generic vanishing. (A similar statement holds for compact K\"ahler manifolds, cf. \cite{pp4}.) The existence of the sheaf $\widehat{\OO_X}$ (or equivalently the perversity of $\R\Phi_P \omega_X$) has significant consequences via the results in \S5, as we will see in \S6.3.

In the reverse direction, the non-vanishing of higher direct images implies a more precise version of Corollary \ref{nef_perverse} in the case of the canonical bundle. 

\begin{corollary}\label{canonical_perverse}
The object $\R \Phi_P \omega_X$ is a perverse sheaf if and only if $X$ is of maximal Albanese dimension.
\end{corollary}
\begin{proof}
The if part follows from Corollary \ref{nef_perverse}. Assume now that $\R\Phi_P \omega_X$ is perverse. We have seen 
that this is equivalent to the fact that $R^i \Phi_{P^\vee} \OO_X = 0$ for $i \neq d$. We claim on the other hand that a result of 
Koll\'ar implies that if the Albanese map $a$ of $X$ has generic fiber of dimension $k$, then $R^{d-k} \Phi_{P^\vee} \OO_X \neq  0$. Indeed, by \cite{kollar} Theorem 3.1 one has a decomposition in the derived category
$$\R a_* \omega_X \cong \oplus_{i =0}^k R^i a_* \omega_X [-i].$$
This decomposition implies in particular that $H^d (X, \omega_X) \cong H^{d-k} (A, R^k a_* \omega_X)$, as the support of 
all $R^i a_* \omega_X$ is at most $d-k$ dimensional. Hence we must have $R^k a_* \omega_X \neq 0$. But then the claim follows 
by a standard application of Grothendieck duality Lemma \ref{gd}.
\end{proof}

Finally, the extra results listed in \S6.1 produce other examples of $WIT_d$-objects, so via Fourier-Mukai duals of perverse
sheaves. For instance:

\noindent
$\bullet$ The same statements as in Corollary \ref{dual} hold replacing $\OO_X$ by $L^{-1}$ for any nef line bundle
$L$ on $X$.

\noindent
$\bullet$ Under the hypotheses of Theorem \ref{direct}, $A_j: = \R\Delta (R^j f_* \omega_X \otimes L)$ satisfies $WIT_d$ for each $j$, so $\widehat{A_j} =  \R\Phi_P A_j [d]$ are all sheaves.

\noindent
{\bf 6.3. Bounding the holomorphic Euler characteristic and applications to irregular varieties.}
The syzygy filtration in \S5, combined with the Evans-Griffith theorem, can be applied to bound the holomorphic Euler characteristic of irregular varieties (i.e. those with $h^1(\OO_X) \neq 0$, or equivalently $b_1(X)\neq 0$) of maximal Albanese dimension. Starting with such a variety $X$, we consider the sheaf $\widehat{\OO_X}$ as in Corollary \ref{dual}, and find the suitable ${\rm Per}_m(Y)$ to which its dual $\R\Phi_P \omega_X$ belongs according to the existence of nontrivial irregular fibrations of 
$X$.

\begin{corollary}\label{holo_chi}
If $X$ is a smooth projective complex variety of maximal Albanese dimension:

\noindent
(i) $\R\Phi_P \omega_X \in {\rm Per}_{gv(\omega_X)} (\Pic0)$. 

\noindent
(ii) $\chi(\omega_X)  \ge gv(\omega_X)$.
\end{corollary}
\begin{proof}
This follows immediately from Theorem \ref{syzygy_fm} and Corollary \ref{inequality} applied to 
$\omega_X$. The only thing that needs mention is that $gv (\omega_X) < \infty$, as 
$0 \in V^d_P (\omega_X)$ with $d = \dim X$.
\end{proof}

The key point is that according to \cite{gl2},  $gv (\omega_X)$ is bounded in terms of the fibrations of $X$ over lower dimensional varieties of maximal Albanese dimension (those that cannot be fibered further should be considered the building blocks in the study of irregular varieties).

\begin{corollary}
Let  $X$ be of maximal Albanese dimension, with $q =  h^1(X, \OO_X)$ and $d = \dim X$. If $X$
is not fibered over any nontrivial normal projective variety of maximal Albanese dimension 
satisfying $q (\tilde Y) - \dim Y \ge q - d  - m + 1$ for any smooth model $\tilde Y$, then 
$$\R\Phi_P \omega_X \in {\rm Per}_m (\Pic0) {\rm~and~} \chi(\omega_X) \ge m.$$
\end{corollary}
\begin{proof}
One applies Corollary \ref{holo_chi}, noting that the hypothesis implies  $gv (\omega_X) \ge m$. 
Indeed, an argument based on \cite{gl2} Theorem 0.1 says that $X$ is fibered over a normal 
projective variety $Y$ of lower dimension such that any smooth model $\tilde Y$ is of maximal 
Albanese dimension, and 
$$q(X) - \dim X- (q(\tilde Y) - \dim  Y)\le gv(\omega_X).$$  
The argument is described in the proof of \cite{pp4} Theorem B.
\end{proof}

The most significant instance of this result is an extension to arbitrary dimension of a classical result on surfaces due to Castelnuovo and de Franchis. I only state a slightly weaker version here for simplicity. Note that it holds in the K\"ahler case as well.

\begin{theorem}[\cite{pp4} Theorem A]\label{cdf} 
Let $X$ be an irregular smooth projective complex variety. If $X$ does not admit any surjective morphism with connected fibers onto a normal projective variety $Y$ with $0 < \dim Y<\dim X$ and with any smooth model $\tilde Y$ of maximal Albanese dimension, then
$$\chi(\omega_X) \ge q(X)-\dim X.$$
\end{theorem}

Along similar lines, arguments based on Theorem \ref{syzygy_fm} lead to classification results. For instance, they are crucial in 
another recent extension to arbitrary dimension due to Barja-Lahoz-Naranjo-Pareschi of a well-known fact on the bicanonical map of surfaces of general type.
 
\begin{theorem}[\cite{blnp} Corollary B]
Under the same assumptions as in Theorem \ref{cdf}, if in addition $X$ is of general type and ${\dim}~X < q(X)$, the linear system 
$|2K_X|$ induces a birational map unless $X$ is birational to a theta divisor in a principally polarized abelian variety.
\end{theorem}

\noindent
{\bf 6.4. Perverse sheaves coming from special subvarieties of abelian varieties.}
Let $(A,\Theta)$ be an indecomposable principally polarized abelian variety (ppav) of dimension $g$, 
and let $P$ be a Poincar\'e bundle on $A \times \widehat{A}$. For subvarieties $X\subset A$,  the question whether $\R\Phi_P (\I_X(\Theta))$ is a perverse sheaf on $\widehat{A}$ is related to 
a beautiful geometric problem. First, note that some of the most widely studied special subvarieties of ppav's satisfy this property:

\noindent
$\bullet$ If $(J(C), \Theta)$ is the polarized Jacobian of a curve of genus $g \ge 2$, and for any $1\le d \le g-1$ we denote by $W_d$ the image of the $d$-th symmetric product of $C$ in $J(C)$ via an Abel-Jacobi map, then  $\I_{W_d} (\Theta)$ is a geometric 
$GV$-sheaf, i.e. $\R\Phi_P (\I_X(\Theta))$ is perverse (cf. \cite{pp1} Proposition 4.4).

\noindent
$\bullet$ If $X \subset \PP^4$ is a smooth cubic hypersurface, $(J(X), \Theta)$ is its intermediate 
Jacobian (a five dimensional ppav which is not the Jacobian of a curve), and $F \subset J(X)$ is the 
Fano surface parametrizing lines on $X$, then $\I_{F} (\Theta)$ is a geometric $GV$-sheaf, i.e. $\R\Phi_P (\I_F(\Theta))$ is perverse (cf. \cite{hoering} Theorem 1.2).

\noindent
Note also that further results along these lines in the case of Prym varieties can be found in  \cite{clv}.

What the $W_d$'s and $F$ have in common is that they are subvarieties whose cohomology classes
are \emph{minimal} (i.e. not divisible in $H^* (A, \ZZ)$). Indeed, $[W_d] = \frac{\theta^{g-d}}{(g-d)!}$ and 
$[F] = \frac{\theta^3}{3!}$ by results of Poincar\'e and Clemens-Griffiths respectively. An problem stemming from work of Beauville and Ran in low dimensions, and suggested in general by Debarre \cite{debarre}, is to show that these are in fact the only subvarieties of ppav's representing minimal cohomology classes. In \cite{pp3} Conjecture A, it is proposed that this should also be equivalent 
to the fact that $\I_X (\Theta)$ is a $GV$-sheaf (so since it is obviously geometric, to the fact that $\R\Phi_P (\I_X(\Theta))$ is a perverse sheaf). One direction is known: in \cite{pp3} Theorem B it is shown that if 
$X$ is a nondegenerate closed reduced subscheme of pure dimension $d$ of a ppav $(A, \Theta)$ of dimension $g$
and  $\I_X (\Theta)$ is a $GV$-sheaf, then $X$ is Cohen-Macaulay and $[X] = \frac{\theta^{g-d}}{(g-d)!}$.
The main ingredient is the criterion provided by Theorem \ref{negative} above.

\noindent
{\bf 6.5. Moduli spaces of vector bundles.}
The most natural higher rank analogue of generic vanishing involving ${\rm Pic}^0(X)$ is to consider the (singular) moduli space
$M_X(r)$ of semistable rank $r$ vector bundles with trivial Chern classes on a smooth projective $X$. 
Some bounds on the dimension of cohomological support loci in $M_X(r)$
are given by Arapura \cite{arapura} \S7, but a more thorough 
understanding remains a very interesting problem. (Note however that the \emph{structure} of these loci is very well understood 
in \cite{arapura}, by means of non-abelian Hodge theory.)

On the other hand, generic vanishing statements can be considered whenever $M$ is a fine moduli space of objects on $X$, 
with a universal object $\E$ on $X\times M$ inducing the functor 
$$\R \Phi _{\E}: \D(X) \longrightarrow \D(M).$$ 
In other words, the problem is to study the variation of the cohomology of the objects parametrized by $M$.
In practice, in order to apply Theorem \ref{negative} the main
difficulty to be overcome is a good understanding of the transform
$\widehat{A^{-1}} = \R \Psi_{\E} (A^{-1})[{\rm dim}~M]$,
with $A$ a very positive line bundle on $M$. Few concrete examples seem to be known: 
besides the case of curves, there are only sparse examples of
generic vanishing phenomena on moduli spaces in higher dimension, coming from constructions of 
Yoshioka on $K3$ and abelian surfaces and of Bridgeland-Maciocia on threefolds with fibrations
by such surfaces. I finish with a very brief discussion of these.

\noindent
{\bf Curves.}
The issue whether $\R\Phi_{\E} \omega_X$ is a perverse sheaf on moduli spaces of higher rank vector bundles on a curve $X$ is not so interesting, as it is easy to check by hand. However, Theorem \ref{negative} is useful when applied to bundles other than 
$\omega_X$. For instance, when applied to higher rank stable bundles, it links the entire indeterminacy loci of determinant line bundles on these moduli spaces to a well-known construction of Raynaud, as first considered by Hein \cite{hein}. This is of a somewhat different flavor from the main direction of this note, so I will refrain from including details. The interested reader can consult \cite{pp2} \S7.

\noindent
{\bf Some surfaces and threefolds.}
Consider first $X$ to be a complex abelian or $K3$ surface.  For a coherent sheaf $E$ on $X$,
the  Mukai vector of $E$ is
$$v(E) : = {\rm rk} (E) + c_1(E) + (\chi (E) -  \epsilon \cdot rk(E) ) [X] \in H^{ev}(X, \ZZ),$$
where $\epsilon$ is $0$ if $X$ is abelian and $1$ if $X$ is $K3$.
Given a polarization $H$ on $X$ and a vector $v \in H^{ev}(X, \ZZ)$,
we consider the moduli space $M_H (v)$ of sheaves $E$ with $v(E)
= v$, stable with respect to $H$. If the Mukai vector $v$ is
primitive and isotropic, and $H$ is general, the moduli space is $M
= M_H (v)$ is smooth, projective and fine, and it is in fact again
an abelian or $K3$ surface (cf. e.g. \cite{yoshioka1}). The
universal object $\E$ on $X \times M$ induces an equivalence of
derived categories $\R \Phi_{\E} : \D(X) \rightarrow \D(M)$.
Yoshioka gives many examples that amount to satisfying:

\noindent ($*$) ~ The Mukai vector $v$ is primitive and isotropic, and the structure
sheaf $\OO_X (\cong \omega_X)$ satisfies $WIT_2$ with respect to $\R \Phi_{\E^\vee}$ (or equivalently 
is $GV$ with respect to  $\R \Phi_{\E}$ by Theorem \ref{negative}).

\noindent
$\bullet$ Let $(X, H)$ be a polarized $K3$ surface such
that ${\rm Pic}(X) = \ZZ\cdot H$, with $H^2 = 2n$. Let $k > 0$ be an
integer such that $kH$ is very ample. Consider $v = k^2n + kH + [X]$. 
It is shown in \cite{yoshioka3} Lemma 2.4 that under these assumptions 
$(*)$ is satisfied.\footnote{It is also shown in \emph{loc. cit.} \S2 that in
fact $X \cong M$.}

\noindent
$\bullet$ Let $(X, H)$ be a polarized abelian
surface with ${\rm Pic}(X) = \ZZ\cdot H$. Write $H^2 = 2r_0 k$, with
$(r_0,k) = 1$. Consider the Mukai vector $v_0 := r_0 + c_1(H) +
k[X]$. By \cite{yoshioka2}, Theorem 2.3 and the
preceding remark, in this situation $(*)$ is again satisfied by $\OO_X$ 
(among many other examples).

A Calabi-Yau fibration is a morphism $\pi:
X\rightarrow S$ of smooth projective varieties, with connected
fibers, such that $K_X\cdot C=0$ for all curves $C$ contained in
fibers of $\pi$. Assuming that $X$ is a threefold, it
is an elliptic, abelian surface, or $K3$-fibration (in the sense
that the nonsingular fibers are of this type). Say $\pi$ is flat,
and consider a polarization $H$ on $X$, and $Y$ an irreducible
component of the relative moduli space $M^{H, P}(X/S)$ of sheaves on
$X$ (over $S$), semistable with respect to $H$, and with fixed
Hilbert polynomial $P$. The choice of $P$ induces on every smooth
fiber $X_s$ invariants which are equivalent to the choice of a Mukai
vector $v \in H^{ev}(X_s, \ZZ)$ as above. Assuming that $Y$ is
also a threefold, and fine, Bridgeland and Maciocia (cf. \cite{bm},
Theorem 1.2) proved that it is smooth, and the induced morphism
$\hat \pi: Y \rightarrow S$ is a Calabi-Yau fibration of the same
type as $\pi$. If $\E$ is a universal sheaf on $X\times
Y$, then $\R \Phi_{\E}: \D(X) \rightarrow
\D(Y)$ is an equivalence. 
Now if the first half of ($*$) is satisfied for each smooth fiber of $\pi$, 
it is proved in \cite{bm} \S7 that the moduli space $M^{H, P}(X/S)$ does have a fine
component $Y$ which is a threefold, so the above applies. The same method 
based on Theorem \ref{negative} leads to the following:

\begin{proposition}[\cite{pp2}, Proposition 7.7]\label{cy}
Let $X$ be a smooth projective threefold with a smooth Calabi-Yau
fibration $\pi: X \rightarrow S$ of relative dimension two. Let $H$
be a polarization on $X$ and $P$ a Hilbert polynomial, and assume
that condition ($*$) is satisfied for each fiber of $\pi$. Consider a fine
three-dimensional moduli space component $Y \subset M^{H,P} (X/S)$,
and let $\E$ be a universal sheaf on $X\times Y$. Then $\omega_X $
is a $GV_{-1}$-sheaf with respect to $\E$. In particular
$$H^i (X, \omega_X\otimes E) = 0, ~{\rm for~all} ~i > 1 ~{\rm and~all}~ E\in Y {\rm~general}.$$
\end{proposition}

This can be extended to the case 
of singular fibers, involving a slightly technical condition which 
should nevertheless be often satisfied; see \cite{pp2} Remark 7.8.
A similar statement holds for threefold elliptic fibrations.

\section{Appendix: some homological commutative algebra}

A useful technical point is that over Cohen-Macaulay rings one can avoid 
finite homological dimension hypotheses in statements of Auslander-Buchsbaum-type, by 
involving the canonical module.  The key fact is the following consequence of Grothendieck 
duality for local cohomology.

\begin{lemma}[\cite{bh}, Corollary 3.5.11]\label{cm_numerical}
Let $(R, m)$ be a local Cohen-Macaulay ring of dimension $n$ with canonical module $\omega_R$, and let $M$ be a finitely 
generated module over $R$ of depth $t$ and dimension $d$. Then:

\noindent
(a) ${\rm Ext}^i_R (M, \omega_R) = 0$ for $i < n-d$ and $i > n-t$.

\noindent
(b) ${\rm Ext}^i_R (M, \omega_R) \neq 0$ for $i = n-d$ and $i = n-t$.

\noindent
(c) ${\rm dim}~{\rm Ext}^i_R (M, \omega_R) \le n-i$ for all $i \ge 0$.
\end{lemma}

\begin{corollary}\label{support}
Let $\F$ be a coherent sheaf on a Cohen-Macaulay scheme of finite type over a regular local ring, 
or a complex analytic space $X$. Then:

\noindent
(a) $\E xt^i (\F , \omega_X) = 0 ~ for~ all~ i < {\rm codim}~{\rm Supp}~\F$ 
and ~$\E xt^i (\F , \omega_X) \neq 0 ~ for~ i = {\rm codim}~{\rm Supp}~\F$.

\noindent
(b) ${\rm codim}~ {\rm Supp} ~ \E xt^i (\F , \omega_X) \ge i$ for all $i$.

\noindent
(c) If $\F$ is locally free, then  $\E xt^i (\F , \omega_X) = 0$ for all $i >0$.
\end{corollary}

Let now $X$ be a noetherian scheme of finite type over a field, or a complex analytic space.

\begin{definition}
Let $\F$ be a coherent sheaf on $X$. Then $\F$ is called a \emph{$k$-th syzygy sheaf} if locally there exists an exact sequence 
\begin{equation}\label{syzygy}
0\longrightarrow \F \longrightarrow \E_k \longrightarrow \ldots \longrightarrow \E_1 \longrightarrow \G \longrightarrow 0
\end{equation}
with $\E_j$ locally free for all $j$. It is well known for example that if $X$ is normal, then $1$-st syzygy sheaf is equivalent to torsion-free, and $2$-nd syzygy sheaf is equivalent to reflexive. Every coherent sheaf is declared to be a $0$-th syzygy sheaf. A locally free one is declared to be an $\infty$-syzygy sheaf.
\end{definition}

\noindent
Following \cite{ab} and \cite{eg1}, we consider a variant of Serre's condition $S_k$.

\begin{definition}\label{serre}
A coherent sheaf $\F$ on $X$ satisfies property $S_k^{\prime}$ if for all $x$ in the support of $\F$ we have:
$${\rm depth}~ \F_x \ge {\rm min} \{k, {\rm dim}~\OO_{X, x}\}.$$
\end{definition}

\begin{proposition}\label{cm_fundamental}
Let $X$ be Cohen-Macaulay, and let $\F$ be a coherent sheaf on $X$. Consider the following conditions:

\noindent
(a) $\F$ is a $k$-th syzygy sheaf.

\noindent
(b) ${\rm codim}~ {\rm Supp} ~  \mathcal{E}xt^i (\F, \omega_X) \ge i + k, {\rm ~for~ all~} i > 0.$

\noindent
(c) $\F$ satisfies $S_k^{\prime}$.

\noindent
Then (b) is equivalent to (c), and if in addition $\F$ has finite homological dimension or $X$ is Gorenstein 
in codimension less than or equal to one, they are both equivalent to (a).
\end{proposition}
\begin{proof}
When $\F$ is of finite homological dimension or $X$ is Gorenstein in codimension less than or equal to one, the equivalence between (a) and (c)  is a well-known result of Auslander-Bridger, \cite{ab} Theorem 4.25.
The equivalence between (b) and (c) is also standard in the case of finite homological dimension, when it follows from the Auslander-Buchsbaum formula (cf.  \cite{hl} Proposition 1.1.6(ii)). We note 
here that a variation of the usual argument, involving Lemma \ref{cm_numerical}, proves it in general.
Consider first $x \in  {\rm Supp} ~  \mathcal{E}xt^i (\F, \omega_X)$ for some $i >0$. Then by Lemma
\ref{cm_numerical}(i) we get that $i \le \dim ~\OO_{X, x} - {\rm depth}~\F_x$. Combined with $S_k^{\prime}$, this implies $\dim \OO_{X, x} \ge i +k$. On the other hand, consider $x \in {\rm Supp}~\F$. We need to show that $\dim \OO_{X, x} - {\rm depth}~\F_x\le {\rm max} \{\dim \OO_{X, x} - k , 0\}$.
But by Lemma \ref{cm_numerical} (i) and (ii), we have 
$$\dim \OO_{X, x} - {\rm depth}~\F_x = {\rm max}~ \{p ~|~ {\rm Ext}^p_{\OO_{X, x}} (\F_x, \omega_{X, x}) \neq 0\}.$$
But for any such $p$ we have $x \in  {\rm Supp} ~  \mathcal{E}xt^p (\F, \omega_X)$. Hence either the 
maximum is $0$ and we're done, or it is positive and for all such $p >0$ we have $\dim \OO_{X, x} \ge p+ k$.
\end{proof}

\noindent
Finally, a key result here is the Syzygy Theorem of Evans-Griffith.

\begin{theorem}[\cite{eg1} Corollary 1.7]\label{syzygy_theorem}
Let $X$ be a Cohen-Macaulay scheme over a field, and $\F$ a $k$-th syzygy sheaf of finite homological dimension on $X$ which is not locally free. Then ${\rm rank}(\F) \ge k$.
\end{theorem}

\providecommand{\bysame}{\leavevmode\hbox
to3em{\hrulefill}\thinspace}

\end{document}